\documentclass[11pt,a4paper]{article}

\usepackage[utf8]{inputenc}
\usepackage[T1]{fontenc}
\usepackage[francais]{babel}
\usepackage{amsmath}
\usepackage{amsfonts}
\usepackage{amssymb}
\usepackage{amsthm}
\usepackage{mathtools}
\usepackage{pgf,tikz}
\usepackage{mathrsfs}
\usepackage{float}
\usetikzlibrary{arrows}
\usepackage{xcolor}
\usepackage[pdfborder={0 0 0}]{hyperref}
\usepackage{stmaryrd}
\usepackage{cancel}
\usepackage{enumitem}
\usepackage{bbold}
\usepackage{fullpage}
\newtheorem{thm}{Théorème}
\newtheorem*{thm2}{Théorème}
\newtheorem{corollaire}[thm]{Corollaire}
\newtheorem{lemma}[thm]{Lemme}

\theoremstyle{definition}
\newtheorem*{remer}{Remerciements}
\newtheorem{rem}{Remarque}

\newcommand{\N}{\mathbb{N}}

\newcommand{\R}{\mathbb{R}}
\newcommand{\C}{\mathbb{C}}

\newcommand{\e}{\mathrm{e}}
\newcommand{\dt}{\mathrm{d}}

\title{Théorème d'Erd\H{o}s-Kac dans un régime de grande déviation pour les translatés d'entiers ayant $k$ facteurs premiers}
\author{Olivier GARCONNET}
\date{}

\newcommand{\Addresses}{ \bigskip
  \footnotesize

  O.~Garçonnet, \textsc{Aix Marseille Univ, CNRS, I2M, Marseille, France}\par\nopagebreak
  \textit{E-mail addresse:} \texttt{olivier.garconnet@univ-amu.fr}

}

\begin{document}

\maketitle

\begin{abstract}
Let $x\geqslant 3$, and for $n$ an integer, let $\omega(n)$ be its number of distinct prime factors. We show that, among the values $n\leqslant x$ with $\omega(n)=k$ where $1\leqslant k \ll \log_2 x$, $\omega(n-1)$ satisfies an Erdős-Kac type theorem around $2\log_2 x$, so in large deviation regime, when weighted by~$2^{\omega(n-1)}$.
This sharpens a result of Gorodetsky and Grimmelt \cite{GG24} with  a quantitative and quasi-optimal error term. 
The proof of the main theorem is based on the characteristic function method and uses recent progress on Titchmarsh's divisor problem, in particular a theorem of Fouvry and Tenenbaum \cite{FT22}. 
\end{abstract}

\section{Présentation des résultats}
Pour $k, n \geqslant 1$ deux entiers et $x \geqslant 1$ un réel, on note $\omega(n)$ le nombre de facteurs premiers distincts de $n$, on note $\mathcal{E}_k(x)$ l'ensemble suivant
$$
\mathcal{E}_k(x):=\{n \leqslant  x: \quad \omega(n)=k\},
$$
et $\pi_k(x):=\#\mathcal{E}_k(x)$. On définit
$$
\Phi(y):=\frac{1}{\sqrt{2 \pi}} \int_{-\infty}^y \mathrm{e}^{-t^2 / 2} \mathrm{d} t . \quad(y \in \mathbb{R})
$$

Le théorème suivant est un résultat fondamental en théorie probabiliste des nombres. Il donne la répartition de la fonction $\omega$ autour de sa valeur moyenne selon une loi normale. 

\begin{thm2}[Erdős \& Kac \cite{EK39}; Rényi \& Turán \cite{RT58}; voir chapitre III.4.4 de \cite{Ten22}] Uniformément pour~$x \geqslant~3$ et $y \in \mathbb{R}$, on a
$$
 \#\left\{n \leqslant  x: \quad \omega(n) \leqslant  \log _2 x+y \sqrt{\log _2 x}\right\}=x\left\{\Phi(y)+O\left(\frac{1}{\sqrt{\log _2 x}}\right)\right\}.
$$ 
\end{thm2}
 
Il est naturel de savoir si un sous-ensemble des entiers défini par une propriété multiplicative, lorsqu'il est translaté, conserve des propriétés de factorisation. 
Halberstam \cite{Hal56} a étudié la répartition de $\omega(p-1)$ pour $p$ un nombre premier et a obtenu un analogue du théorème d'Erdős-Kac.
Fouvry et Tenenbaum~\cite{FT96} prouvent un résultat analogue pour les entiers friables. Goudout~\cite{Gou20}, en s'inspirant de leur méthode, a étudié les entiers avec un nombre fixé de facteurs premiers. 

\begin{thm2}[Goudout \cite{Gou20}]
Soit $R>0$ fixé. Uniformément pour $x \geqslant 3,1 \leqslant  k \leqslant  R \log _2 x$ et~$y \in \mathbb{R}$, on a
$$
\pi_k(x, y)=\sum_{n\in\mathcal{E}_k(x)} \mathbb{1}_{\omega(n-1) \leqslant  \log _2 x+y \sqrt{\log _2 x}}=\pi_k(x)\left\{\Phi(y)+O\left(\frac{\log _3 x}{\sqrt{\log _2 x}}\right)\right\} .
$$
\end{thm2}

Estimer $\omega(n-1)$ pour $n\in\mathcal{E}_k(x)$, loin de la moyenne, $\omega(n-1)\approx\beta\log_2 x$ avec $\beta\neq1,$ est un problème plus compliqué, comme l'indique un théorème d'Elliott \cite[theorem~6]{Ell15}. Cela revient essentiellement à savoir étudier $\sum \beta^{\omega(n-1)}$. 
Le cas~$\beta=1$ revient au théorème de Goudout~\cite{Gou20}. On s'attend à avoir des résultats pour $\beta=2$ en vertu du problème de Titchmarsh comme étudié par Fouvry \cite{Fou85}.
C'est l'objet de la conjecture B de \cite{Ell15} prouvée par Gorodetsky et Grimmelt~\cite{GG24}.

\begin{thm2}[Gorodetsky \& Grimmelt \cite{GG24}]
Soit $y\in\R$. Alors on a
$$\frac{1}{\varphi(N)}\sum_{\substack{p<N \\ \omega(N-p) \leqslant 2 \log \log N+y(2 \log \log N)^{1 / 2}}} 2^{\omega(N-p)} \rightarrow \frac{1}{\sqrt{2 \pi}} \int_{-\infty}^y \e^{-u^2 / 2} \dt u,$$
pour $N\rightarrow\infty$, avec $\varphi$ l'indicatrice d'Euler et la somme porte sur les nombres premiers.
\end{thm2}

Un résultat similaire est attendu, avec les mêmes méthodes, en remplaçant $N-p$ par $p-1$ dans l'énoncé de ce théorème.
Pour prouver ce résultat, les auteurs de \cite{GG24} se basent sur la méthode des moments.

Le but de cet article est de mettre à profit les méthodes analytiques qui sous-tendent \cite{Gou20} afin d'obtenir un terme d'erreur quantitatif et de généraliser le résultat aux entiers $n$ ayant un nombre fixé de facteurs premiers.
On définit

$$
\begin{aligned}
S_k(x):=&\sum_{n\in\mathcal{E}_k(x)} 2^{\omega(n-1)},\\
S_k(x,y):=&\sum_{n\in\mathcal{E}_k(x)} 2^{\omega(n-1)}\mathbb{1}_{\omega(n-1) \leqslant 2\log _2 x+y \sqrt{2\log _2 x}}.
\end{aligned}$$
Le théorème principal de ce papier est le suivant.

\begin{thm}\label{th1}
Soit $R>0$ fixé. Uniformément pour $x\geqslant 3$, $y\in\R$ et $1\leqslant k\leqslant R\log_2x$ on a 
$$S_k(x,y)=S_k(x)\left\{\Phi(y)+O\left(\frac{\log_3x}{\sqrt{\log_2x}}\right) \right\}.$$
\end{thm}
L'estimation de $S_k(x)$, qui est une variante du problème des diviseurs de Titchmarsh pour~$\mathcal E_k$~\cite{Khr98}, sera obtenue au cours de la démonstration et fournira 
$$S_k(x)=\frac{x\left(\log _2 x\right)^{k-1}}{(k-1) !}\left\{A_r+O\left(\frac{1}{\log_2x}\right)\right\},$$
avec $r=\frac{k-1}{\log_2 x}$ et
$$A_r:=\frac{1}{\Gamma( r+1)} \prod_{p \geqslant 2}\left(1+\frac{r+1}{p-1}\right)\left(1-\frac{1}{p}\right)^{ r+1 }.$$

Notre méthode, avec quelques ajustements techniques, permet de préciser le résultat de~\cite{GG24} pour~$N-p$ à la place de $p-1$, et permet d'obtenir un terme d'erreur optimal au facteur~$\log_3 N$~près.
\begin{corollaire}\label{cor1}
Uniformément pour $N\geq 3$ et $y\in\R$, on a
$$\frac{1}{\varphi(N)}\sum_{\substack{p<N \\ \omega(N-p) \leqslant 2 \log_2 N+y(2 \log_2 N)^{1 / 2}}} 2^{\omega(N-p)} =\Phi(y)+O\left(\frac{\log _3 N}{\sqrt{\log _2 N}}\right),$$
où la somme porte sur les nombres premiers.
\end{corollaire}

Pour démontrer ces résultats, on étudie la série génératrice adéquate à l'aide d'un résultat de Fouvry et Tenenbaum \cite{FT22}. Puis on utilise la méthode des fonctions caractéristiques et l'inégalité de Berry-Esseen. De plus, en corollaire de la preuve on démontre le résultat suivant qui localise le nombre de ``petits'' facteurs premiers de~$n-1$ pour~$n\in\mathcal{E}_k(x)$.

%On démontre ce théorème, après avoir étudié la série génératrice adéquate, par la méthode des fonctions caractéristiques et l'inégalité de Berry-Esseen. De plus, en corollaire de la preuve on démontre le résultat suivant qui localise le nombre de ``petits'' facteurs premiers de~$n-1$ pour~$n\in\mathcal{E}_k(x)$.
Pour $k\geqslant 1,$ $\ell\geqslant 0$ et $1\leqslant w\leqslant x,$ on note 
$$
\begin{aligned}
\omega(n, w)  &:=\sum_{\substack{p \mid n \\ p \leqslant  w}} 1, \quad(n \geqslant 1),\\
S_{k,\ell}(x,w)&:=\sum_{n\in\mathcal{E}_k(x)} 2^{\omega(n-1)}\mathbb{1}_{\omega(n-1,w)=\ell}.
\end{aligned} $$

\begin{thm}\label{th2}
Il existe une constante absolue $\eta>0$ telle que pour $R\geqslant 1$ fixé quelconque, uniformément pour $3 \leqslant  w \leqslant  x^{\eta(\log_3x)/(R^2\log_2x)}$, $1 \leqslant  k \leqslant  R \log _2 x$ et $0 \leqslant  \ell \leqslant  R \log _2 w$, on a
$$
S_{k, \ell}(x, w)=S_k(x) \frac{\left(2\log _2 w\right)^{\ell}}{\ell ! (\log w)^2}\left\{H_k\left(\frac{\ell}{\log _2 w}\right)+O\left(\frac{k}{\left(\log _2 x\right)^2}+\frac{\ell+1}{\left(\log _2 w\right)^2}\right)\right\}.
$$
où $H_k$ est la fonction entière définie sur $\mathbb{C}$ par 
$$H_k(z):=\e^{\gamma(2z-2)}\prod_{p}\left(1+\frac{2z-2}{p+r}\right)\left(1-\frac{1}{p}\right)^{2z-2}.$$
et $\gamma$ est la constante d'Euler. 
\end{thm}
\begin{rem}\label{moment}
La méthode utilisée permet de calculer les moments relatifs au Théorème \ref{th1}. Pour tout entier $m\geqslant 0$ fixé, on a 
\begin{align*}
    \sum_{n \in \mathcal{E}_k(x)} 2^{\omega(n-1)}\left(\frac{\omega(n-1)-2 \log _2 x}{\sqrt{2 \log _2 x}}\right)^m ={} \left\{
    \begin{array}{ll}\displaystyle
         O\left(S_k(x)\frac{\log_3 x}{\sqrt{\log_2 x}}\right),& \text{ si $m$ est impair,}  \\
         \displaystyle S_k(x)\left(\frac{(2i)!}{2^i i!}+O\left(\frac{\log_3 x}{\sqrt{\log_2 x}}\right)\right),& \text{ si $m=2i$ est pair.}
    \end{array}\right.
\end{align*}
On retrouve asymptotiquement les moments de la loi normale centrée réduite.

\end{rem}
\subsubsection*{Organisation du papier}

La section \ref{sec etude serie gen} est consacrée à l'étude de la série génératrice associée à la localisation des «~petits~» facteurs premiers de $n-1$. On s'appuie sur un résultat de Fouvry et Tenenbaum~\cite{FT22}. Dans la section \ref{preuve th}, on utilise cette étude et le théorème de Berry-Esseen pour prouver nos résultats.

\section{Étude de la série génératrice}{\label{sec etude serie gen}}
Pour $1\leqslant w\leqslant x$ et $k\geqslant 1$, on définit 
$$F_{k}(z)=F_{k}(z,x,w):=\sum_{n\in\mathcal{E}_k(x)}2^{\omega(n-1)}z^{\omega(n-1, w)}.$$
On omet les variables $w$ et $x$ afin d'alléger la notation. La fonction $F_k(z)$ est la série génératrice associée à~$\left(S_{k,\ell}(x,w)\right)_\ell$ et est calculée dans le lemme suivant.

\begin{lemma}{\label{lem serie generatrice}}
Il existe $\delta >0$ une constante absolue telle que pour $R\geqslant 1$ fixé, $3\leqslant w \leqslant  x,|z| \leqslant R$,  et~$k\leqslant R\log_2 x$, avec $r:=(k-1)/\log_2 x$ et $u:=\log x/\log w$, on ait
$$
\begin{aligned}
F_k(z)=&\;S_{k}(x)(\log w)^{2 z-2}\left\{H_k(z)+\frac{r\xi(z)}{2\log_2 x}+O\left(\frac{1}{\log w}+\frac{\log w}{\log x}+\frac{k^2}{(\log_2 x)^4}\right)\right\}
\\ &+ O\left(xu^{-\delta u}\left(\log x\right)^{2R+4}+x^{1-\delta/4}\right),
\end{aligned}
$$
avec $\xi$ une fonction entière pouvant dépendre des paramètres $w,x$ et $k$, admettant $1$ comme zéro et uniformément bornée pour $|z|\leqslant R$.
\end{lemma}

\begin{rem}
Dans le cas particulier où $k=1$, donc pour les nombres premiers, on a un terme d'erreur plus précis  en $O\left(\frac{1}{\log x}\right)$ au lieu de $O\left(\frac{r^2}{(\log_2 x)^2}\right)$ mais cela n'améliore pas le terme d'erreur du Théorème \ref{th1}.
\end{rem}
Dans le reste de cette section, on démontre le Lemme \ref{lem serie generatrice}.

\subsubsection*{Introduction de $g_z$}

On introduit $g_z$ la fonction multiplicative définie pour un nombre premier $p$ par le tableau suivant :

$$
\begin{array}{|c|c|c|}
\hline
 &  p \leqslant w & p > w\\
\hline
g_z(p)  & 2(z - 1)& 0 \\
\hline
g_z(p^2)  & 1 - 2z& -1 \\
\hline
g_z(p^\alpha) \text{ pour } \alpha > 2 & 0 & 0\\
\hline
\end{array}$$

Cette fonction satisfait $g_z * \tau(n)=2^{\omega(n)}z^{\omega(n,w)}$ pour tout $n\in\N$.
Ainsi on a 

\begin{equation}\label{eq:F_k}
F_k(z)=\sum_{n\in\mathcal{E}_k(x)}2^{\omega(n-1)}z^{\omega(n-1, w)}=\sum_{n\in\mathcal{E}_k(x)}
\sum_{\substack{q\mid n-1}} g_z(q) \tau\left(\frac{n-1}{q}\right).    
\end{equation}

\subsubsection*{Découpage pour $q$ grand}
On montre que l'on peut tronquer la somme à $q\leqslant x^\delta$, où $\delta>0$ est un paramètre arbitraire choisi plus tard. On pose
\begin{equation}\label{eq:Fk1}
F_k(z)=F_k^{(1)}(z) +R_1,\quad\text{où}\quad R_1:=\sum_{n\in\mathcal{E}_k(x)} \sum_{\substack{q|n-1 \\ q>x^{\delta}}} g_z(q)\tau\left(\frac{n-1}{q}\right).
\end{equation}\newline
{}
Par définition de $g_z$, on écrit de façon unique $n-1=q\times q_3=q_1\times q_2^2 \times q_3$, avec $P^+(q_1)\leqslant w$ et~$P^-(q_2)>w$, et on a 
$$
\begin{aligned}
|R_1|\leq\sum_{n\in\mathcal{E}_k(x)} \sum_{\substack{q|n-1 \\ q>x^{\delta}}} |g_z(q)|\tau\left(\frac{n-1}{q}\right) &\leqslant   \sum_{\substack{q_1q_2^2>x^{\delta}\\P^+(q_1)\leqslant w\\ P^-(q_2)>w\\q_1q_2^2q_3<x}} (2R+2)^{\omega(q_1)}\tau\left(q_3\right)\\
&  \ll x\log x\sum_{\substack{q_1q_2^2>x^{\delta}\\P^+(q_1)\leqslant w\\ P^-(q_2)>w}} \frac{(2R+2)^{\omega(q_1)}}{q_1q_2^2}.
\end{aligned}
$$
Avec l'astuce de Rankin, pour $\alpha\in[0,1/3]$, on a 
$$R_1\ll x\log x\sum_{\substack{P^+(q_1)\leqslant w\\ P^-(q_2)>w}} \left(\frac{q_1q_2^2}{x^\delta}\right)^\alpha\frac{(2R+2)^{\omega(q_1)}}{q_1q_2^2}= \frac{x\log x}{x^{\delta\alpha}}\prod_{p\leqslant  w}\left(1+\frac{2R+2}{p^{1-\alpha}-1}\right)\prod_{p>w}\left(1+\frac{1}{p^{2(1-\alpha)}-1}\right).$$\newline
{}
Le deuxième produit converge absolument. Pour le premier produit, le lemme $2$ de \cite{Ten84} implique l'inégalité suivante :
$\sum_{p\leqslant w}\frac{1}{p^{1-\alpha}}\ll \log\log w+w^\alpha$ donc avec $u:=\frac{\log x}{\log w}$ et $\alpha=\min\left(\log (2u)/\log w, 1/3\right)$, on obtient 
\begin{equation}\label{eq:R1}
\begin{aligned}
R_1&\ll x\log x\exp\Big(-\delta u(\log 2u)+(2R+2)2u+(2R+2)\log\log w\Big)+x^{1-\delta/4}\\&\ll xu^{-\delta u/2}\left(\log x\right)^{2R+3}+x^{1-\delta/4}.
\end{aligned}
\end{equation}

\subsubsection*{Préliminaire en vue de l'application du théorème de Fouvry-Tenenbaum}{\label{FT}}
On a l'égalité suivante 
$$\tau(n)= 2\sum_{\substack{d\mid n\\d<\sqrt{n}}}1+\mathbb{1}_{n \text{ carré}}.$$
En accord avec cette égalité on écrit $F_k^{(1)}(z)=F_k^{(2)}(z)+R_2$ où
\begin{equation}\label{eq:R2}
R_2:=\sum_{n\in\mathcal{E}_k(x)}\sum_{\substack{q\mid n-1\\q\leqslant x^\delta}} g_z(q)\mathbb{1}_{\frac{n-1}{q} \text{ carré}}\ll x^\varepsilon\sum_{\substack{q,\ell\\q\leqslant x^\delta\\\ell^2q+1\leqslant x}}1\ll x^{\varepsilon+\delta+1/2}.    
\end{equation}
La majoration utilise le fait que $|g_z(q)|\ll_\varepsilon x^\varepsilon$. On doit à présent estimer 
 \begin{equation}\label{eq:Fk2}
F_k^{(2)}(z)=\sum_{\substack{q \leqslant x^\delta }} 2g_z(q) \sum_{\substack{d\\d^2 q+1\leqslant x}} \sum_{\substack{d^2 q + 1 \leqslant n \leqslant x \\ n \equiv 1 [dq]}} \mathbb{1}_{\omega(n)=k}.  
\end{equation}

On souhaite appliquer le théorème $1.5$ de Fouvry-Tenenbaum \cite{FT22} pour estimer la somme en~$n$ et enlever la condition de congruence sur $n$. Pour cela on a besoin d'une dépendance multiplicative en $n$, donc d'après une idée de Selberg \cite{Sel54} on introduit 
$$\mathbb{1}_{\omega(n)=k}=\frac{1}{2\pi i}\int_{|v|=1} \frac{1}{v^{k+1}} v^{\omega(n)} \dt v.$$
\newline
{}
On décompose $F_k^{(2)}(z)=F_k^{(3)}(z)+R_3$, où
\begin{align}
\label{eq:R3} F_k^{(3)}(z)={}&\sum_{\substack{q \leqslant x^\delta }} 2g_z(q) \sum_{\substack{d\\d^2 q+1\leqslant x}}\frac{1}{\varphi(dq)}\sum_{\substack{n\in\mathcal{E}_k(x)\\ d^2 q + 1 \leqslant n \\ (n,dq)=1}} 1,\\
R_3:={}&\frac{1}{2\pi i}\int_{|v|=1} \frac{1}{v^{k+1}} \sum_{\substack{q \leqslant x^\delta }} 2g_z(q) \sum_{\substack{d\\d^2 q+1\leqslant x}} \left(\sum_{\substack{d^2 q + 1 \leqslant n \leqslant x \\ n \equiv 1 [dq]}} v^{\omega(n)}-\frac{1}{\varphi(dq)}\sum_{\substack{d^2 q + 1 \leqslant n \leqslant x \\ (n,dq)=1}} v^{\omega(n)}\right)\dt v.\notag\end{align}
\newline
{}
Les résultats de Fouvry-Tenenbaum \cite{FT22} nécessitent l'indépendance des bornes des sommes. On écrit donc $R_3=R_3'-R_3''$ où dans $R_3^\prime$ la somme sur $n$ porte sur les entiers $n\leqslant x$ et dans $R_3^{\prime\prime}$ elle porte sur les entiers $n\leqslant d^2q+1$.

\subsubsection*{Séparation des variables et application du théorème de Fouvry-Tenenbaum}
On détaille la majoration de $R_3^{\prime}$; la même méthode s'appliquera pour $R_3^{\prime\prime}$. Pour supprimer la dépendance en $q$ et $d$ des bornes, on utilise un principe de séparation des variables : on se place dans des intervalles $d\in [\lambda^a,\lambda^{a+1}[$ et $q\in [\lambda^b,\lambda^{b+1}[$ avec~$\lambda:=1+(\log x)^{-B}$, et $B>1$. On écrit d'abord

$$R_3^{\prime}\ll\sup_{|v|=1}|R_3^{\prime}(v)|,$$
où
$$R_3^{\prime}(v):=\sum_{q\leqslant x^\delta}\sum_{\substack{ d^2 q +1\leqslant x}}g_z(q)\left(\sum_{\substack{n \leqslant x \\ n \equiv 1[dq]}} v^{\omega(n)}-\frac{1}{\varphi(q d)} \sum_{\substack{n \leqslant x \\(n, d q)=1}} v^{\omega(n)}\right).$$\newline
{}
On choisit $N$ tel que $\lambda^N \leqslant x^\delta < \lambda^{N+1}$. On trouve alors que
\begin{equation}\label{eq:R3primeprime}
\begin{aligned}
\left|R_3^{\prime}(v)\right|\leq & \left|\sum_{\substack{0\leqslant b\leqslant N,\\0\leqslant a,\, \lambda^{2a+b+3} \leqslant x}} \sum_{\substack{\lambda^b \leqslant q<\min(\lambda^{b+1},x^\delta)\\ \lambda^a \leqslant d<\lambda^{a+1}}}g_z(q)\Biggl(\sum_{\substack{n \leqslant x \\ n=1[d q]}} v^{\omega(n)}-\frac{1}{\varphi(q d)} \sum_{\substack{n \leqslant x \\(n, d q)=1}} v^{\omega(n)}\Biggr)\right|
\\& +\sum_{\substack{0\leqslant b\leqslant N, \,0\leqslant a\\ \lambda^{2a+b} \leqslant x<\lambda^{2a+b+3} }} \sum_{\substack{\lambda^b\leqslant q<\min(\lambda^{b+1},x^\delta)\\ \lambda^a \leqslant d<\lambda^{a+1}}}|g_z(q)|\left|\sum_{\substack{ n\leqslant x\\ n=1[d q]}} v^{\omega(n)}-\frac{1}{\varphi(q d)} \sum_{\substack{ n\leqslant x  \\(n, d q)=1}} v^{\omega(n)}\right|.
\end{aligned}
\end{equation}\newline
{}
 On choisit $\delta<1/105$, puis on utilise le théorème de Fouvry-Tenenbaum~\cite{FT22} sur le premier terme noté (PT) de \eqref{eq:R3primeprime}. On obtient une majoration de ce terme en $$PT\ll_A\sum_{\substack{0\leqslant b\leqslant N \\0\leqslant a\\ \lambda^{2a+b} \leqslant x}}\frac{x}{(\log x)^A}.$$\newline
{}
On somme sur $a$ puis sur $b$ et on obtient une majoration en $$PT\ll \frac{Nx}{(\lambda-1)(\log x)^{A-1}}\ll \frac{x}{(\log x)^{A-2B-2}}.$$

Pour le second terme noté (ST) de \eqref{eq:R3primeprime}, on majore trivialement ce qui se trouve dans les valeurs absolues par $\ll \frac{x}{dq}$, on obtient

$$ST\ll x\sum_{\substack{0\leqslant b\leqslant N,\,0\leqslant a\\ \lambda^{2a+b} \leqslant x<\lambda^{2a+b+3}}} \sum_{\substack{\lambda^b\leqslant q<\min(\lambda^{b+1},x^\delta)\\ \lambda^a \leqslant d<\lambda^{a+1}}}\frac{|g_z(q)|}{dq}.$$\newline
{}
On somme ensuite sur $d$ puis sur $a$, il reste
$$ST \ll x(\lambda-1)\sum_{0\leqslant b\leqslant N } \sum_{\lambda^b\leqslant q<\min(\lambda^{b+1},x^\delta)}\frac{|g_z(q)|}{q}=x(\lambda-1)\sum_{0<q\leqslant x^\delta}\frac{|g_z(q)|}{q}.$$
Cette dernière somme peut être majorée par $\exp\left(\sum_{p\leqslant x^\delta}\frac{2R+2}{p}\right)$, donc par la formule de Mertens on obtient une majoration du second terme de \eqref{eq:R3primeprime} en $ST\ll x(\log x)^{-B+2R+2}$.\newline
{}
Finalement si on choisit $A$ et $B$ correctement on obtient pour $C>0$ quelconque
$R_3^{\prime}\ll \frac{x}{(\log x)^C}.$

Pour la majoration de $R_3^{\prime\prime}$, on utilise de la même manière un principe de séparation des variables pour pouvoir appliquer le corollaire $1.7$ de Fouvry-Tenenbaum~\cite{FT22}.
On commence par majorer trivialement les éléments pour lesquels $n$ est plus petit que $x^{1/2}$; pour les autres, la contrainte que $n\leqslant qd^2+1$ implique que $\log (qd^2+1)\asymp\log x$. Quitte à supposer $\delta<1/210 $, on obtient ainsi\newline
{}
$$R_3^{\prime\prime}:=\frac{1}{2\pi i}\int_{|v|=1} \frac{1}{v^{k+1}} \sum_{\substack{q \leqslant x^\delta }} 2g_z(q) \sum_{\substack{ d\\d^2q<x}} \left(\sum_{\substack{n \leqslant d^2q \\ n \equiv 1 [dq]}} v^{\omega(n)}-\frac{1}{\varphi(dq)}\sum_{\substack{n \leqslant d^2q \\ (n,dq)=1}} v^{\omega(n)}\right)\dt v\ll_C \frac{x}{(\log x)^C}.$$
où $C$ est une constante arbitraire.

Suite à la majoration de $R_3^\prime$ et celle de $R_3^{\prime\prime}$, on obtient pour $C>0$ quelconque
\begin{equation}\label{eq:R3majo}
R_3\ll_C \frac{x}{(\log x)^C}.    
\end{equation}

\subsubsection*{Bilan intermédiaire et introduction de $f_z$}

Dans l'équation \eqref{eq:R3}, on intègre sur la variable $v$ pour se ramener à étudier les nombres avec $k$ facteurs premiers distincts. De la même manière que la majoration de $R_1$ en \eqref{eq:R1}, avec l'astuce de Rankin, on peut relâcher la condition $q\leqslant x^\delta$ :

\begin{equation}\label{eq:R4}
R_4:= \sum_{\substack{q > x^\delta}} 2g_z(q) \sum_d \sum_{\substack{d^2 q + 1 \leqslant n \leqslant x \\ (n, q d)=1\\ n\in \mathcal{E}_k(x)}} \frac{1}{\varphi(q d)}\ll xu^{-\delta u/2}\left(\log x\right)^{2R+3}+x^{1-\delta/4}.    
\end{equation}\newline
{}
D'après \eqref{eq:Fk1}, \eqref{eq:Fk2} et \eqref{eq:R3}, il nous reste 

\begin{equation}\label{eq:FkFk5}
\begin{aligned}
F_{k}(z)=\sum_{n\in\mathcal{E}_k(x)}2^{\omega(p-1)}z^{\omega(p-1, w)} &=\sum_{\substack{q}} 2g_z(q) \sum_d \frac{1}{\varphi(q d)} \sum_{\substack{d^2 q + 1 \leqslant n \leqslant x \\ (n, q d)=1}} \mathbb{1}_{\omega(n)=k}+R_1+R_2+R_3-R_4
\\&=2\sum_{\ell\leqslant x} f_z(\ell) \sum_{n\in\mathcal{E}_k(x)} \mathbb{1}_{(\ell,n)=1}-R_5+R_1+R_2+R_3-R_4,
\end{aligned}\end{equation}\newline
{}
où $R_5$ et la fonction multiplicative $f_z$ sont définis par

$$
f_z(\ell):=\sum_{\substack{d, q\\d^2 q=\ell}} \frac{g_z(q)}{\varphi(d q)}, \qquad
R_5:=2\sum_{\ell\leqslant x} f_z(\ell) \sum_{\substack{n \leqslant \ell \\ \omega(n)=k}} \mathbb{1}_{(\ell,n)=1}.
$$\newline
{}
Pour la fin de la preuve du Lemme \ref{lem serie generatrice}, on note le terme principal 
\begin{equation}\label{eq:Fk5}
F_k^{(5)}(z)=2\sum_{\ell\leqslant x} f_z(\ell) \sum_{n\in\mathcal{E}_k(x)} \mathbb{1}_{(\ell,n)=1}.   \end{equation}

\subsubsection*{Moyenne d'une fonction multiplicative sur \texorpdfstring{$\mathcal E_k(x)$}{Ek(x)} pour estimer $F_k^{(5)}(z)$}

Le lemme 3 de \cite{Gou20} appliqué à la fonction $n \mapsto\mathbb{1}_{(\ell,n)=1}$ donne, avec $r=\frac{k-1}{\log_2 x}$ (on rappelle que $k\ll \log_2 x$),

\begin{equation}\label{eq:moyfctmult}
\sum_{n \in \mathcal{E}_k(x)} \mathbb{1}_{(\ell,n)=1}=\frac{x}{\log x} \frac{\left(\log _2 x\right)^{k-1}}{(k-1) !}\left\{\lambda_{\ell}(r)-\frac{r \lambda_{\ell}^{\prime \prime}(r)}{2 \log _2 x}+O\left(\frac{k^2}{\left( \log _2 x\right)^4}\right)\right\},
\end{equation}\newline
{}
où $\lambda_{\ell}$ est la fonction entière définie par
$$
\lambda_{\ell}(y):=\frac{1}{\Gamma( y+1)} \prod_{p \geqslant 2}\left(1+\frac{y}{p-1}\right)\left(1-\frac{1}{p}\right)^{ y } \prod_{p|\ell}\left(1+\frac{y}{p-1}\right)^{-1} \quad(y \in \mathbb{C}),$$
où les pôles du second produit en $1-p$ sont compensés par les zéros du premier.
On pose 
$$C_r:=\frac{1}{\Gamma( r+1)} \prod_{p \geqslant 2}\left(1+\frac{r}{p-1}\right)\left(1-\frac{1}{p}\right)^{ r },$$\newline
{}
et dans l'optique d'évaluer $F_k^{(5)}(z)$ à l'aide de \eqref{eq:moyfctmult}, on veut estimer 
$$T_{x}(z,r):=2C_r\sum_{\ell\leqslant x}f_z(\ell)\prod_{p|\ell}\left(1+\frac{r}{p-1}\right)^{-1},$$\newline
{}
ainsi que
$$\sum_{\ell\leqslant x}f_z(\ell)\frac{r\lambda_\ell^{\prime\prime}(r)}{2\log_2x}=\frac{r}{2\log_2x}\frac{\dt^2 T_x(z,r)}{\dt r^2}.$$

\subsubsection*{Application d'une formule de Perron effective pour évaluer $T_x(z,r)$}
Evalué en une puissance de $p$, $f_z$ vaut pour $\alpha\geq0$
$$\begin{aligned}
f_z(p^{2 \alpha+1})= &\begin{cases}\frac{2z-2}{p^\alpha(p-1)} & \text { si } p \leqslant  w \\ 0 & \text { si } p>w,\end{cases}\\
f_z(p^{2 \alpha+2})= &\begin{cases}\frac{1}{p^{\alpha}(p-1)}+\frac{1-2z}{p^{\alpha+1}(p-1)} & \text { si } p \leqslant  w \\ \frac{1}{p^{\alpha}(p-1)}-\frac{1}{p^{\alpha+1}(p-1)} & \text { si } p>w.\end{cases}
\end{aligned}$$\newline
{}
Comme $g:\ell\mapsto \ell f_z(\ell)\prod_{p|\ell}\left(1+\frac{r}{p-1}\right)^{-1}$ est une fonction multiplicative, on introduit son produit eulérien qui vaut

$$\begin{aligned}  &\prod_{p \leqslant w}\left(1+\frac{1}{p^{2 s-1}-1}\frac{p-1}{p-1+r}+\frac{2z-2}{(p-1+r) p^{s-1}}+\frac{(2z-2)\left(1-p^{s-1}\right)}{p^{s-1}(p-1+r)\left(p^{2 s-1}-1\right)}\right)
\\& \times\prod_{p>w}\left(1+\frac{1}{p^{2 s-1}-1}\frac{p-1}{p-1+r}\right) \\ & 
=:\zeta(2s-1)H(s), \end{aligned}$$
où $H(s)$ est bien défini pour $\Re(s)>1/2$ car 
$$H(s)=\prod_{p \leqslant w}\left(1+\frac{2z-2}{p^{s}}+O\left(\frac{1}{p^{2s}}+\frac{1}{p^{3s-1}}\right)\right) \prod_{p>w}\left(1-\frac{1}{p^{2s-1}}\frac{r}{p-1+r}\right).$$\newline
{}
On applique la première formule de Perron effective \cite[théorème II.2.3]{Ten22} à la fonction $g$ et cela nous donne pour $T\geqslant 1$ et $\kappa > 1$ choisi plus tard
\begin{equation}\label{eq:perron}
\sum_{\ell\leqslant x}g(\ell)=\frac{1}{2 \pi i} \int_{\kappa-i T}^{\kappa+i T} \zeta(2s-1)H(s)x^s \frac{\mathrm{d} s}{s}+O\left(x^\kappa \sum_{\ell \geqslant 1} \frac{\left|g(\ell)\right|}{\ell^\kappa(1+T|\log (x / \ell)|)}\right).
\end{equation}
\textbf{Théorème des résidus et majoration de l'intégrale :}
Le rectangle $\big(\kappa-i T;\kappa + i T;\sigma+i T;\sigma-i T\big)$, avec $\kappa=1+\frac{1}{\log x}$ et $\sigma=1-\frac{c}{\log(2+T)}$, se situe dans la région sans zéro de la fonction~$\zeta$. 
De la même manière que pour \eqref{eq:R1}, on a pour $s=1-\varepsilon +i\tau$, $\varepsilon\in [-1/\log x; c/\log(2+T)]$ et~$\tau \in\R,$
$$ |H(s)x^s|\ll  x\exp\left(\sum_{p\leqslant w}\frac{2R+2}{p^{1-\varepsilon}}-\varepsilon\log x\right)\ll  x(\log w)^{2R+2}\exp \left((2R+2)w^\varepsilon -\varepsilon\log x\right).$$\newline
{}
La fonction en $\varepsilon$ atteint son maximum aux bornes de l'intervalle.
Le théorème des résidus autour du contour précédemment cité donne 
\begin{align*}
\frac{1}{2 \pi i} &\int_{\kappa-i T}^{\kappa+i T} \zeta(2s-1)H(s) x^s \frac{\mathrm{d} s}{s}=\frac{xH(1)}{2}\\&+O\left(x(\log w)^{2R+2}\left(\frac{\log T}{T}+(\log T)^2\exp\Big((2R+2)\exp\Big(\frac{c\log w}{\log T}\Big)-\frac{c\log x}{\log T}\Big)\right)   \right).   
\end{align*}\newline
{}
\textbf{Majoration du terme d'erreur de la formule de Perron :} Si on écrit $\ell=ab^2$ avec $a$ sans facteur carré, on remarque que l'on a $|g(\ell)|\leqslant (2R+2)^{\omega(ab)}b$.\newline
Dans la somme du terme d'erreur, la contribution des entiers $\ell$ n'appartenant pas à \newline $I:=[(1-1/\sqrt{T})x,(1+1/\sqrt{T})x]$ est majorée par
$$\frac{x}{\sqrt{T}} \sum_{\ell \notin I} \frac{\left|g(\ell)\right|}{\ell^\kappa}\ll \frac{x}{\sqrt{T}}\sum_{a,b}\frac{(2R+2)^{\omega(ab)}}{(ab)^{\kappa}}\ll \frac{x}{\sqrt{T}}\frac{1}{(\kappa-1)^{2(2R+2)}}\ll \frac{x(\log x)^{2(2R+2)}}{\sqrt{T}}.$$
Pour la contribution des entiers $l\in I$ de la somme du terme d'erreur, on applique tout d'abord l'inégalité de Cauchy-Schwarz pour supprimer les contraintes liées à la fonction $g$. Puis on applique le principe de l'hyperbole, on obtient pour tout $y>0$, la majoration suivante 
$$
\begin{aligned}
x\sum_{\ell\in I} \frac{\left|g(\ell)\right|}{\ell^\kappa}\ll &\, x \sum_{ab^2\in I}\frac{(2R+2)^{\omega(ab)}}{ab}\ll x(\log x)^C\left(\sum_{\substack{a<y\\ab^2\in I}}\frac{1}{ab}+\sum_{\substack{b<\sqrt{x/y}\\ab^2\in I}}\frac{1}{ab}\right)^{1/2}\\ 
\ll &\, x(\log x)^{C}\left(\frac{\log x}{\sqrt{T}}+\frac{\sqrt{y}}{\sqrt{x}}+\frac{1}{y}\right)^{1/2}\ll \frac{x(\log x)^{C+1/2}}{T^{1/4}}+x^{11/12},    
\end{aligned}$$
en choisissant $y=x^{1/3}$ et pour une constante $C=C(R)>0$ suffisamment grande. 

\noindent
\textbf{Choix adéquat de T :} Soit $A:=\exp\left(\sqrt{\log x\log_2 x}\right)$, si $w\leqslant A$ on prend $T=A$,  si $w>A$ on prend $T=\exp\left(\frac{c\log w}{\log (2u)}\right)$, avec $u=\frac{\log x}{\log w}$. Alors il existe $\delta>0$ telle que on a 

$$\sum_{\ell\leqslant x}g(\ell)=\frac{xH(1)}{2}+O\left(x(\log x)^{2R+4}\left(\exp(-\delta u\log u)+\exp\left(-c\sqrt{\frac{\log x}{\log_2x}}\right)\right)+x^{11/12}\right).$$\newline
{}
\textbf{Calcul de $H(1)$ :} On a $$H(1)=\prod_{p\leqslant w}\left(1+\frac{2z-1}{p-1+r}\right)\left(1-\frac{1}{p}\right)\left(1+O\left(\frac{1}{\log w}\right)\right).$$\newline
{}
D'après la troisième formule de Mertens, uniformément pour $ |z|\leqslant R$, on a 

$$\left(\e^\gamma \log w\right)^{2z-2} \prod_{p \leqslant  w}\left(1-\frac{1}{p}\right)^{2z-2}=1+O\left(\frac{1}{\log w}\right),$$
ainsi
$$H(1)=\left( \log w\right)^{2z-2}\e^{\gamma(2z-2)}\prod_{p}\left(1+\frac{2z-1}{p-1+r}\right) \left(1-\frac{1}{p}\right)^{2z-1}\left(1+O\left(\frac{1}{\log w}\right)\right).$$\newline
{}
\textbf{Conclusion :} Finalement on obtient\newline
{}
\begin{equation}\label{eq:sumnfz}
\begin{aligned}
\sum_{\ell\leqslant x}\ell f_z(\ell)\prod_{p|\ell}\left(\frac{p-1}{p-1+r}\right)&=\sum_{l\leqslant x} g(l)\\&=\frac{x}{2}(\log w)^{2z-2}\left(h_k(z)+O\left(\frac{1}{\log w}+\exp\left(-c\sqrt{\frac{\log x}{\log_2 x}}\right)\right)\right)\\
&\hspace{13em}+O\left(x(\log x)^{2R+4}u^{-\delta u}+x^{11/12}\right),
\end{aligned}
\end{equation}
avec $c>0$ une constante, $u=\frac{\log x}{\log w}$ et 
$$h_r(z):=\e^{\gamma(2z-2)}\prod_{p}\left(1+\frac{2z-1}{p-1+r}\right)\left(1-\frac{1}{p}\right)^{2z-1}.$$ \newline
{}
Puis par intégration par partie, en faisant attention à la condition $w\leqslant x$, on obtient
\begin{equation*}
\begin{aligned}
\sum_{n\leqslant x}f_z(n)\prod_{p|n}\left(\frac{p-1}{p-1+r}\right)=
\frac{\log x}{2} (\log w)^{2z-2}&\left(h_r(z)+O\left(\frac{1}{\log w}+\frac{\log w}{\log x}+\exp\left(-c\sqrt{\frac{\log x}{\log_2 x}}\right)\right)\right) \\
&\qquad\qquad\qquad\quad+O\left(x(\log x)^{2R+4}u^{-\delta u}+x^{11/12}\right).
\end{aligned}
\end{equation*}
Finalement on a uniformément en $r\leqslant R$
$$
\begin{aligned}
T_x(z,r)=C_r(\log x)(\log w)^{2z-2}&\left(h_r(z)+O\left(\frac{1}{\log w}+\frac{\log w}{\log x}+\exp\left(-c\sqrt{\frac{\log x}{\log_2 x}}\right)\right)\right)\\
&\qquad\qquad\qquad\qquad\qquad+O\left(x(\log x)^{2R+4}u^{-\delta u}+x^{11/12}\right).
\end{aligned}$$

\subsubsection*{Conclusion de la preuve du Lemme \ref{lem serie generatrice}}

Dans la partie précédente, l'équation \eqref{eq:sumnfz} avec l'estimation de Goudout \eqref{eq:moyfctmult} permet de calculer $R_5$,

\begin{align}
R_5&=2\sum_{\ell\leqslant x} f_z(\ell) \sum_{\substack{n \leqslant \ell \notag\\ \omega(n)=k}} \mathbb{1}_{(\ell,n)=1}\\&= \frac{C_r}{\log x} \frac{\left(\log _2 x\right)^{k-1}}{(k-1) !}\sum_{\ell\leqslant x} \ell f_z(\ell)\prod_{p|\ell}\left(\frac{p-1}{p-1+r}\right)\left\{1+O\left(\frac{1}{\log_2x}\right)\right\} + O_\varepsilon(x^\varepsilon)\notag\\
&=  \frac{1}{\log x} \frac{C_r x\left(\log _2 x\right)^{k-1}}{(k-1) !}(\log w)^{2z-2}\left\{h_r(z)+O\left(\frac{1}{\log w}+\frac{1}{\log_2x}\right)\right\}+O\left(x(\log x)^{2R+4}u^{-\delta u}+x^{11/12}\right).\label{eq:R5}
\end{align}

Pour conclure on obtient
\begin{equation}\label{eq:Fk5esti}
\begin{aligned}
F_k^{(5)}(z)=\frac{x\left(\log _2 x\right)^{k-1}}{(k-1) !}(\log w)^{2z-2}&\left\{C_rh_r(z)-\frac{r}{2\log_2x}\frac{\dt^2 C_rh_r(z)}{\dt r^2}+O\left(\frac{1}{\log w}+\frac{\log w}{\log x}+\frac{k^2}{(\log_2x)^4}\right)\right\}\\&\qquad\qquad\qquad\qquad\qquad\qquad+O\left(x(\log x)^{2R+4}u^{-\delta u}+x^{11/12}\right).
\end{aligned}    
\end{equation}

Finalement on regroupe les estimations \eqref{eq:FkFk5}, \eqref{eq:Fk5esti}, \eqref{eq:R1}, \eqref{eq:R2}, \eqref{eq:R3majo}, \eqref{eq:R4} et \eqref{eq:R5}. En évaluant l'équation précédente en $z=1$ et $w=\exp\big((\log_2 x)^4\big)$, on obtient l'approximation suivante de~$S_k(x)=F_k(1)$, où on a posé $A_r=C_rh_r(1)$,

$$S_k(x)=\frac{x\left(\log _2 x\right)^{k-1}}{(k-1) !}\left\{A_r-\frac{r}{2\log_2x}\frac{\dt^2 C_rh_r(1)}{\dt r^2}+O\left(\frac{k^2}{(\log_2x)^4}\right)\right\}.$$\newline
{}
Ce qui permet d'avoir l'écriture exacte du terme souhaité
$$
\begin{aligned}
F_k(z)=S_{k}(x)(\log w)^{2z-2}&\left\{\frac{h_r(z)}{h_r(1)}+\frac{r}{2\log_2 x}\left(\frac{h_r(z)\frac{\dt^2 C_rh_r(1)}{\dt r^2}}{C_rh_r(1)^2}-\frac{\frac{\dt^2 C_rh_r(z)}{\dt r^2}}{C_rh_r(1)}\right)\right. \\&\qquad+\left. O\left(\frac{1}{\log w}+\frac{\log w}{\log x}+\frac{k^2}{(\log_2 x)^4}\right)\right\}+O\left(x(\log x)^{2R+4}u^{-\delta u}+x^{1-\delta/4}\right).
\end{aligned}$$
On pose pour tout $z\in\C$

$$
\begin{aligned}
H_k(z):=&\frac{h_k(z)}{h_k(1)}=\e^{\gamma(2z-2)}\prod_{p}\left(1+\frac{2z-2}{p+r}\right)\left(1-\frac{1}{p}\right)^{2z-2},\\
\xi(z):=&\frac{h_r(z)}{C_rh_r(1)^2}\frac{\dt^2}{\dt r^2}\Bigl[C_rh_r(1)\Bigr]-\frac{1}{C_rh_r(1)}\frac{\dt^2 }{\dt r^2}\Bigl[C_rh_r(z)\Bigr].
\end{aligned}$$
La fonction $\xi$ est une fonction entière admettant $1$ pour zéro et bornée uniformément en $w,x$ et~$k$ sous les conditions de l'énoncé du Lemme \ref{lem serie generatrice} et cela conclut la preuve du Lemme \ref{lem serie generatrice}.

\section{Preuve des Théorèmes \ref{th1} et \ref{th2} et preuve du Corollaire \ref{cor1}}\label{preuve th}
\subsection{Preuve du Théorème \ref{th2}}
\begin{proof}
On utilise une idée de Selberg  \cite{Sel54}. La quantité $S_{k,\ell}(x,w)$ est le coefficient de~$z^\ell$ dans $F_k(z)$, ce qui est d'après le théorème des résidus
$$S_{k, \ell}(x, w)=\frac{1}{2 i \pi} \oint_{|z|=\rho} \frac{F_k(z,x,w)}{z^{\ell+1}} \mathrm{~d} z.$$
On insère l'estimation du Lemme \ref{lem serie generatrice} dans l'intégrale en choisissant $\rho=\ell/\log_2 w$. L'intégrale résultante en $z$ est alors estimée par une méthode similaire à la section $3$ de \cite{Gou20}, la seule différence étant un facteur $2$ dans la puissance de $\log w$.
\end{proof}

\subsection{Preuve du Théorème \ref{th1}}\label{preuve th1}
\subsubsection*{Approximation de $S_k(x,y)$ par une version sans les grands facteurs premiers}
On suit la démonstration de \cite{Gou20}, elle-même basée sur la démonstration de \cite[cor.5]{FT96}. Pour cette démonstration, on choisit

$$w:=\exp\left(\frac{\log x}{\left(\log_2 x\right)^2}\right).$$
On se donne une constante $C>0$ à fixer plus tard. On pose pour $x\geqslant 2$ et $y\in \R$

$$\begin{aligned} \tilde{S}_k(x,y) & :=\sum_{n \in \mathcal{E}_k(x) }2^{\omega(n-1)}\mathbb{1}_{\omega(n-1, w) \leqslant  2\log _2 x+y \sqrt{2\log _2 x}}, \\
D_k(x) & :=\sum_{n \in \mathcal{E}_k(x) }2^{\omega(n-1)}\mathbb{1}_{\omega(n-1)-\omega(n-1, w)>C \log _3 x}.
\end{aligned}$$
Pour montrer que l'influence des grands facteurs premiers de $n-1$ est négligeable, donc que~$\tilde{S}_k(x,y)$ est une bonne approximation de $S_k(x,y)$, il suffit de majorer $D_k$, car on a 

\begin{equation}\label{eq:encadrement Sk}
\tilde{S}_k \left(x,y-\frac{C \log _3 x}{\sqrt{\log _2 x}}\right) -D_k(x)\leqslant  S_k(x, y) \leqslant  \tilde{S}_k\left(x,y\right). 
\end{equation}\newline
{}
On utilise l'astuce de Rankin pour majorer $D_k$, ce qui se traduit par l'inégalité

$$\begin{aligned}D_k(x)\leqslant {} & \sum_{n\in\mathcal{E}_k(x)} 2^{\omega(n-1,w)} 4^{\omega(n-1)-\omega(n-1,w)}2^{-C \log_3 x}\\ \leqslant  {} & \frac{1}{(\log_2 x)^{C\log2}}\sum_{n\in\mathcal{E}_k(x)}2^{\omega(n-1,w)}\tau\bigg(\prod_{\substack{q^\nu \| n-1 \\ q>w}} q^\nu\bigg)^2.\end{aligned}$$
D'après \cite{Lan89}, pour tous $\delta>0$ et $n\geqslant 1$, on a 

$$\tau(n)^2 \ll_\delta \sum_{\substack{q \mid n \\ q \leqslant  n^\delta}} \tau(q)^{2 / \delta}.$$
De plus on a 

$$2^{\omega(n-1,w)}\leqslant  \tau(n-1,w)\leqslant  2\sum_{\substack{d \mid n-1 \\ d \leqslant  \sqrt{n-1}\\P^+(d)\leqslant w}}1,$$ 
ce qui fournit
$$D_k(x)   \ll\frac{1}{(\log_2 x)^{C\log2}} \sum_{\substack{ d \leqslant  \sqrt{x}\\P^+(d)\leqslant w}}\sum_{\substack{ q \leqslant  x^{\delta}\\P^-(q)>w}}\tau(q)^{2/\delta} \sum_{\substack{n\in\mathcal{E}_k(x)\\n\equiv1[dq]}}1. $$

Comme dans la section \ref{sec etude serie gen}, on applique le théorème $1.5$ de Fouvry-Tenenbaum \cite{FT22} afin de supprimer la condition de congruence, quitte à choisir $\delta>0$ suffisamment petit.
On obtient pour toute constante $C^\prime>0$ 
$$R_6:=\frac{2}{(\log_2 x)^{C\log2}} \sum_{\substack{ d \leqslant  \sqrt{x}\\P^+(d)<w}}\sum_{\substack{ q \leqslant  x^{\delta}\\P^-(q)>w}}\tau(q)^{2/\delta}\left|\sum_{\substack{n\in\mathcal{E}_k(x)\\n\equiv1[d]\\n\equiv1[q]}}1-\frac{1}{\varphi(dq)}\sum_{\substack{n\in\mathcal{E}_k(x)\\(n,dq)=1}}1\right|\ll \frac{x}{(\log x)^{C^\prime}}.$$\newline
{}
En utilisant \eqref{eq:moyfctmult}, on déduit

\begin{align}
D_k(x)&\ll\frac{\pi_k(x)}{(\log_2 x)^{C\log2}} \sum_{\substack{ d \leqslant  \sqrt{x}\\P^+(d)\leqslant w}}\sum_{\substack{ q \leqslant  x^{\delta}\\P^-(q)>w}}\tau(q)^{2/\delta} \frac{1}{\varphi(qd)}+R_6 \notag
\\&\ll  \frac{\pi_k(x)}{(\log_2 x)^{C\log2}}\sum_{\substack{ d \leqslant  \sqrt{x}\\P^+(d)\leqslant w}}\frac{1}{\varphi(d)}\left(\frac{\log x}{\log w}\right)^{2^{2/\delta}}+R_6 \notag
\\&\ll \frac{S_k(x)}{\log_2 x}, \label{eq:D_k}
\end{align}\newline
{}
en choisissant $C$ suffisamment grand et car $(\log x) \pi_k(x)\asymp S_k(x)$.

\subsubsection*{Approximation de $\tilde{S}_k(x,y)$ grâce à l'inégalité de Berry-Esseen}
Enfin, pour estimer $\tilde{S}_k(x,y)$, on utilise la méthode des fonctions caractéristiques. Il est naturel de poser 
$$
\tilde{S}_k^{\star}(x, y):=\sum_{n \in \mathcal{E}_k(x) }2^{\omega(n-1)}\mathbb{1}_{\omega(n-1, w) \leqslant  2\log _2 w+y \sqrt{2\log _2 w}}.
$$\newline
{}
Il suffit d'établir la majoration
\begin{equation}\label{eq:berry essen}
\sup _{y \in \mathbb{R}}\left|\tilde{S}_k^{\star}(x, y)-S_k(x) \Phi(y)\right| \ll \frac{S_k(x)}{\sqrt{\log _2 x}}.    
\end{equation}\newline
{}
En effet, on a~$\log_2w=\log_2x-2\log_3x$ et $\tilde{S}_k(x, y)=\tilde{S}_k^{\star}(x, y^\prime)$ avec $y^\prime=y+O\left(\frac{\log_3x}{\sqrt{\log_2 x}}\right)$. De plus,
le résultat \eqref{eq:berry essen} implique que la fonction $y\mapsto\tilde{S}_k^{\star}(x, y)$ est approchée par une gaussienne, qui varie peu avec $y$, plus précisément, on a $||\Phi '||_\infty\ll1$. 
Donc en supposant~\eqref{eq:berry essen} vérifiée, on obtient 
$$\sup _{y \in \mathbb{R}}\left|\tilde{S}_k(x, y)-S_k(x) \Phi(y)\right| \ll \frac{\log_3x}{\sqrt{\log _2 x}}S_k(x).$$\newline
{}
En procédant de même, grâce aux inégalités~\eqref{eq:encadrement Sk}, \eqref{eq:D_k}, on obtient 
$$\sup _{y \in \mathbb{R}}\left|S_k(x, y)-S_k(x) \Phi(y)\right| \ll \frac{\log_3x}{\sqrt{\log _2 x}}S_k(x).$$
Ce qui conclut la preuve du Théorème \ref{th1}.\newline
{}

On revient à la preuve de \eqref{eq:berry essen}. D'après le théorème de Berry-Esseen  \cite[th.II.7.16]{Ten22}, en posant~$T:=2\log _2 w$, on a 
$$\sup _y\left|\tilde{S}_k^{\star}(x, y)-S_k(x)\Phi(y)\right| \ll \frac{S_k(x)}{\sqrt{T}}+\int_{-\sqrt{T}}^{\sqrt{T}}\left|\e^{-i t \sqrt{T}} F_k(\e^{it/ \sqrt{T}})-S_k(x)\e^{-t^2/2}\right| \frac{\dt t}{|t|}.$$\newline
{}
Pour $\vartheta \in \mathbb{R}$, on a $\mathrm{e}^{i \vartheta}=1+O(|\vartheta|)$. Ainsi, lorsque $|t| \leqslant  1 / \log x$,
$$
\begin{aligned}
\left|\e^{-i t \sqrt{T}} F_k\left(\e^{i t / \sqrt{T}}\right)- S_k(x)\e^{-t^2/2}\right| & \ll |t| \sum_{n \in \mathcal{E}_k(x)}2^{\omega(n-1)}\left|\frac{\omega(n-1, w)}{\sqrt{T}}-\sqrt{T}\right|+|t|^2S_k(x) \\
& \ll|t|S_k(x)\left(\frac{\log w}{\sqrt{T}}+\sqrt{T}+|t|\right),
\end{aligned}
$$\newline
{}
où l'on a majoré $\omega(n-1, w) \ll \log w$. On a donc
$$
\int_{-1 / \log x}^{1 / \log x}\left|\e^{-i t \sqrt{T}} F_k\left(\e^{i t / \sqrt{T}}\right)- S_k(x)\e^{-t^2/2} \right| \frac{\dt t}{|t|} \ll \frac{S_k(x)}{\sqrt{T}}\left(\frac{\log w}{\log x}+\frac{T}{\log x}+\frac{\sqrt{T}}{(\log x)^2}\right)\ll \frac{S_k(x)}{\sqrt{T}} .
$$\newline
{}
Pour $1 / \log x \leqslant |t| \leqslant  T^{1 / 6}$, on applique le Lemme \ref{lem serie generatrice}. Puisque $H_k(1)=1$ et $\xi(1)=0$, on déduit

$$
F_k\left(\e^{i t / \sqrt{T}}\right)=S_k(x) \e^{it \sqrt{T}} \e^{-t^2/2 }\left\{1+O\left(\frac{|t|+|t|^3}{\sqrt{T}}+\frac{1}{\left(\log_2x\right)^2}\right)\right\},
$$
et donc
$$
\int_{1 / \log x<|t| \leqslant  T^{1/6}}\left|\e^{-i t \sqrt{T}} F_k\left(\e^{i t / \sqrt{T}}\right)-S_k(x)\e^{-t^2/2 }\right| \frac{\dt t}{|t|} \ll \frac{S_k(x)}{\sqrt{T}}\left(1+\frac{\sqrt{T}}{(\log_2 x)^2}\right) \ll \frac{S_k(x)}{\sqrt{T}}.
$$\newline
{}
Finalement, pour $|t|>T^{1/6}$, on utilise l'inégalité  $\cos (\theta)-1 \leqslant -2 \theta^2 / \pi^2$, valable pour $-1 \leqslant  \theta \leqslant  1$. Ainsi pour $T^{1 / 6}<|t| \leqslant  \sqrt{T}$, le Lemme \ref{lem serie generatrice} fournit la majoration

$$
\left|F_k\left(\e^{i t / \sqrt{T}}\right)\right| \ll S_k(x) \e^{-2 t^2 / \pi^2} .
$$
On obtient donc

$$
\int_{T^{1/6}<|t| \leqslant  \sqrt{T}}\left|\e^{-i t \sqrt{T}} F_k\left(\e^{i t / \sqrt{T}}\right)-S_k(x)\e^{-t^2/2 }\right| \frac{\dt t}{|t|} \ll  S_k(x)\log T\exp(-2T^{1/3}/\pi^2)\ll  \frac{S_k(x)}{\sqrt{T}}.
$$

\subsection{Preuve du Corollaire \ref{cor1}}
\begin{proof}
On étudie cette fois-ci la série génératrice $$F(z)=F(z,N,w):=\sum_{n\in \mathcal{E}_k(N)}2^{\omega(N-n)}z^{\omega(N-n, w)}.$$
On ne peut pas appliquer le théorème de Fouvry-Tenenbaum \cite{FT22} car nous comptons à présent les $n\in \mathcal{E}_k(N)$ et satisfaisant $n\equiv N [qd] $. La classe de congruence $N$ croît trop rapidement avec la taille typique de $n$. Cependant dans le cas des nombres premiers, on peut appliquer à la place la proposition $5.5$ de \cite{GG24} qui provient de \cite[prop. 6.1]{BGLRM23}. De la même manière qu'à la section \ref{sec etude serie gen}, on obtient pour $3\leqslant w \leqslant  N,|z| \leqslant R$ et~$u:=\log N/\log w$,
$$
F(z)=S(N)(\log w)^{2 z-2}\left\{1+O\left(\frac{1}{\log w}+\frac{1}{(\log_2 N)^4}\right)\right\} + O\left(Nu^{-\delta u}\left(\log N\right)^{2R+4}+N^{1-\delta/4}\right),
$$
avec $S(N)=\sum_{p<N}2^{\omega(N-p)}=\varphi(N)\left(1+O\left(\frac{1}{\log_2 N}\right)\right)$.

Ensuite le même travail que dans la section \ref{preuve th1} et le théorème de Berry-Esseen permettent de conclure.
\end{proof}

\begin{rem}
Pour calculer les moments évoqués à la remarque \ref{moment}, on fait apparaître la fonction caractéristique à l'aide d'une intégrale. Pour tout $m\in\N$, on a
\begin{align*}
    \sum_{n \in E_k(x)} 2^{\omega(n-1)}&\left(\frac{\omega(n-1)-2 \log _2 x}{\sqrt{2 \log _2 x}}\right)^m
    \\&= \frac{m!}{2\pi i}\left(\frac{1}{\sqrt{2 \log _2 x}}\right)^m \oint_{\Gamma} \sum_{n \in E_k(x)} 2^{\omega(n-1)} z^{\omega(n-1)-2 \log_2(x)} \frac{\dt z}{z(\log z)^{m+1}},
\end{align*}
où $\Gamma$ est un contour autour de $1$.
On montre que l’influence des grands facteurs premiers de $n-1$ est négligeable de la même manière que dans l'encadrement \ref{eq:encadrement Sk}. Ensuite on applique notre Lemme \ref{lem serie generatrice} et on trouve les moments annoncés dans la remarque \ref{moment}.
\end{rem}

\begin{remer}
Je tiens à exprimer toute ma gratitude à mon directeur de stage, Sary Drappeau, pour son soutien précieux, ses nombreuses relectures et tous ses conseils, qui ont été essentiels à la rédaction de cet article.   
\end{remer}

\bibliographystyle{amsalpha}
\bibliography{bib}
\Addresses

\end{document}